\documentstyle[12pt]{amsart}

\newtheorem{theorem}{Theorem}

\newtheorem{lemma}{Lemma}

\newtheorem{definition}{Definition}
\newtheorem{corollary}{Corollary}
\newtheorem{example}{Example}

\newcommand{\btree}{$\beta(1,0)$-tree }
\newcommand{\btrees}{$\beta(1,0)$-trees }
\newenvironment{proof}{\noindent {\bf Proof:}}{{\qed}}

\newcommand{\ril}{right-to-left }

\newcommand{\vanish}[1]{}

\begin{document}

\title
[2-stack sortable permutations with a given number of runs]
{2-stack sortable permutations\\with a given number of runs}

\author{Mikl\'os B\'ona}
\thanks{This paper was written while the author was
a one-term visitor at  Mathematical Sciences Research Institute in Spring 
1997. This visit
was supported by an MIT Applied Mathematics Fellowship.} 
\address{Department of Mathematics \\
        Massachusetts Institute of Technology \\
        Cambridge, MA 02139}
            
\begin{abstract} Using earlier results we prove a formula for the number
$W_{(n,k)}$ of 2-stack sortable 
permutations of length $n$ with $k$ runs, or in other words, $k-1$
descents. This formula will yield the suprising fact that there are as
many 2-stack sortable permutations with $k-1$ descents as with $k-1$
ascents. We also prove that $W_{(n,k)}$ is unimodal in $k$, for any fixed
$n$.  \end{abstract}

\maketitle

\section{Introduction}
\subsection{Our main results}
In this paper we are going to  show that the number
of 2-stack sortable permutations of length $n$
 with $k-1$ ascents is equal to the number $T_{(n,k)}$ of
\btrees on $n+1$ nodes with $k$ leaves.
 (See Section 2 for 
the definition of $\beta(1,0)$-trees).
 This, and  results from
 \cite{cori} and \cite{jackson} will enable us to easily show that
\begin{equation} \label{formula}
W_{(n,k)}=\frac{(n+k-1)!(2n-k)!}{k!(n+1-k)!(2k-1)!(2n-2k+1)!},\end{equation}
which formula is symmetric in $k$ and $n+1-k$.
\subsection{Background and Definitions}
In what follows, permutations of length $n$ will be called
$n$-permutations. We say that $i$ is a {\em descent} of the permutation
$p=(p_1,p_2,\cdots ,p_n)$ if $p_i>p_{i+1}$. Similarly, $j$ is an ascent
of $p$ if $p_j<p_{j+1}$. If $p$ has $d$ descents, then $p$ decomposes
into $d+1$ increasing sequences of consecutive entries, and we then say
that $p$ has $d+1$ {\em runs}. We will use the concept of descents and
runs interchangeably, according to the current context.

The stack-sorting operation $\Pi$ can be defined on the set of all 
$n$-permutations as follows. Let $p=p_Lnp_R$ be an $n$-permutation, with $p_L$
 and $p_R$ respectively denoting its subword before and after the maximal 
entry. Let $\Pi(p)=\Pi(p_L)\Pi(p_R)n$, where $p_L$ and $p_R$ are defined
 recursively by this same rule. For a nonrecursive, algorithmic definition,
or the origin of the notion see \cite{knuth}, \cite{West}.

A permutation $p$ is called $t$-stack sortable if $\Pi^t(p)$ is the identity
permutation. 

The set of 1-stack sortable $n$-permutations is easy to characterize by the
following notion of pattern avoidance.
Let $q=(q_1, q_2, \cdots,
q_k)$ be a $k$-permutation 
and let $p=(p_1, p_2, \ldots, p_n)$ be an $n$-permutation. We say that
$p$ contains a $q$-subsequence (or $q$-pattern) if there exists $1\leq
i_{q_1}< i_{q_2} < 
\cdots < i_{q_k} \leq n$ such that $p_{i_1}<p_{i_2}<\cdots < p_{i_k}$. We say
 that $p$ avoids $q$ if $p$  contains no $q$-subsequence.

For example, $p$ avoids 231 if it cannot be written as $\cdots,a,\cdots,
b,\cdots c,\cdots $ so that $c<a<b$. It is easy to show \cite{knuth}
 \cite{West} \cite{west2} that a permutation is 1-stack sortable if and
only if it avoids the pattern 231. In paticular, the number of 1-stack sortable
permutations is therefore $C_n$, the $n$-th Catalan number. 

The set of 2-stack sortable permutations is much more complex. For example,
it is not true that a subword of a 2-stack sortable permutation is always
2-stack sortable; for 35241 is 2-stack sortable, while 3241 is not.
Similarly, it is far more difficult to find a formula for the number $W_n$
of  these permutations than for that of 1-stack
sortable ones. See \cite{doron} \cite{goulden} or \cite{dulucq} for a
proof of the formula $W_n=2(3n)!/((n+1)!(2n+1)!)$ and see \cite{West}
for a characterization of 2-stack sortable permutations.

In this paper we are going to construct a bijection between the set of
2-stack sortable $n$-permutations with $k$ runs and that of \btrees on
$n+1$ vertices having $k$ leaves. This bijection will be based on a
bijection in \cite{goulden}, which mapped these permutations into a
class of labeled binary trees. (Our trees are not binary). The set of
these trees 
has recently been shown to be equinumerous \cite{cori} to a certain
class of planar 
maps, and the number of those planar maps is known \cite{jackson}.
Therefore, we are going to obtain a {\em formula} for the number
$W_{(n,k)}$ of 2-stack sortable permutations with $k$ runs, which is a bit
surprising, given that even the formula for $W_n$ was fairly difficult
to prove.

 What is even more surprising is that this formula reveals that
{\em there are as many 2-stack sortable permutations with $k$ descents
as with $k$ ascents}. We do not see any obvious reason why this should
be true. A direct combinatorial proof of this fact would be highly
desirable, but it may not be easy as we do not even have such a proof
for the same fact in  the much simpler case of 132-avoiding
permutations. We also prove that $W_{(n,k)}$ is unimodal in $k$ for any
fixed $n$. 

We point out that this is the second case recently \cite{terkep}
 when a class of the
description trees introduced in \cite{cori} to enumerate planar maps has
been used to enumerate permutations. 

\section{The correspondence between trees and permutations}
Let $p$ be an $n$-permutation.
A {\em \ril maximum} of $p$ is an entry which is larger than all entries it
precedes. If the \ril maxima of $p$ are $a_1>a_2>\cdots >a_t$, then
clearly $a_1=n$ and $a_t$ is the last entry of $p$. In the next three
paragraphs we closely follow the characterization of 2-stack sortable
permutation given in
\cite{goulden}. 

Write $p$ in the form $s_1a_1s_2a_2\cdots
s_ta_t$. Here the $s_i$ are strings of entries located between two \ril
maxima.  We say that $p$ is of type 1 if the entry $a_t-1$ is part of
the string $s_t$. Otherwise we say that $p$ is of type 2. The
1-permutation 1, denoted $\epsilon$, is of type 2.  Then the
following holds.

\begin{lemma} \cite{goulden} \label{egyes}
 There is a bijection $F$ from the set of
2-stack sortable $n$-permutations of type 1 to that of 2-stack sortable
$(n-1)$-permutations with a \ril maximum marked. Moreover, $F$ preserves
the number of descents and does not decrease the number of \ril maxima. 
\end{lemma}
\begin{proof}\cite{goulden} Let $p$ be of type 1, with the above
$s_1a_1s_2a_2\cdots s_ta_t$
decomposition. Define $F$ as follows. Delete  $a_t$ from the end of $p$
and decrement all entries larger than $a_t$ by 1. Then $a_t-1\in s_t$
becomes a \ril maximum in the new permutation, in fact, it becomes the
$t$th \ril maximum. Mark this entry to get $F(p)$. This map is
reversible by simply adding an entry one larger than the marked vertex
to the end of $F(p)$ and increment all the larger entries by one.

 So $F$ is a bijection. Moreover, $F$ does not create or destroy any
existing descents as it deletes the last entry of $p$, and that entry
was larger 
than its predecessor. The number of \ril maxima is not decreased by $F$,
because $a_t-1$ becomes a \ril maximum instead of $a_t$.

Finally, one sees easily \cite{goulden} that $F(p)$ is 2-stack sortable
if and only if $p$ is, and the proof is complete. \end{proof}

We need one more lemma from \cite{goulden}. We omit its proof as we will
not need any of its rather complicated machinery. Let $rl(p)$ be the
number of \ril maxima of $p$ and let $desc(p)$ be the number of descents
in $p$. 

\begin{lemma} \label{decomp} Let $n\geq 2$. Then there is a bijection
$G$ from the set 
of 2-stack sortable $n$-permutations of type 2 onto the set of ordered
pairs $(p_1,p_2)$ where \begin{enumerate}
\item $p_1$ is any 2-stack sortable $n_1$-permutation, $p_2$ is either a
2-stack 
sortable $n_2$-permutation of type 1, or $\epsilon$,
\item  $n_1+n_2=n$,
\item $rl((p_1)+rl(p_2)=rl(p)$,
\item $desc(p_1)+desc(p_2)+1=desc(p)$.
\end{enumerate} \end{lemma} 

We are going to iterate the decomposition of the above Lemma. In other
words, we take a 2-stack sortable $n$-permutation $p$ of type 2,
decompose it into 
$(p_1,p_2)$ by the bijection $G$ of \cite{goulden}, then, if $p_1$ was
of type 2 (and not $\epsilon$), then we apply $G$ to $p_1$ to get
$(p_{11},p_{12})$. Again, if $p_{11}$ 
is of type 2 (and not $\epsilon$), the we apply $G$ to $p_{11}$, and so
on. We stop when both elements of the current decomposition are either
of type 1, or $\epsilon$. As each step of this algorithm uses a
bijection, $p$ can be recovered from its final decomposition. This
proves the following Corollary.  

\begin{corollary} \label{szet} There is a bijection $H$ from the set of
2-stack sortable $n$-permutations onto the set of string $(q_1,q_2,\cdots
,q_r)$, where $r\geq 1$ and \begin{itemize} \item
for all $i$, $q_i$ is either a 2-stack sortable
permutation of type 1, or $\epsilon$ and 
\item $\sum_{i=1}^r|q_i|=n$,
\item $\sum_{i=1}^r rl(q_i)=rl(p)$,
\item the total number of runs in the $q_i$ equals the number of runs in
$p$.
\end{itemize} \end{corollary}
The last part follows from the last part of Lemma \ref{decomp} and the
fact that  the number of runs in any permutation is one more than that
of descents. 

Cori, Jacquard and Schaeffer \cite{cori} give the following definition in
their study of planar maps.
\begin{definition} A rooted plane tree with nonnegative
integer labels $l(v)$ on each of its vertices $v$ is called a
\btree if it satisfies the 
following conditions: 
\begin{itemize} 
\item if $v$ is a leaf, then $l(v)=1$,
\item if $v$ is the root and $v_1,v_2,\cdots, v_k$ are its children,
then $l(v)=\sum_{i=1}^kl(v_k)$,
\item if $v$ is a nonroot
 internal node and $v_1,v_2,\cdots, v_k$ are its children,
then $l(v)\leq  \sum_{i=1}^kl(v_k)$.
\end{itemize}
\end{definition}

\begin{example}{\em Two \btrees are shown in Figures 1a and 1b.}
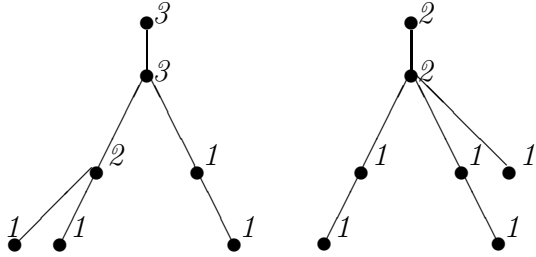
 \[\begin{picture}(230,130)(30,30)

\put(27,54){\line(1,2){33}} \put(93,54){\line(-1,2){33}}
\put(27,54){\circle*{5}} \put(31,57){1} \put(10,54){\line(1,1){29}}
\put(41,81){\circle*{5}} \put(45,84){2} \put(6,57){1} \put(10,54){\circle*{5}}
\put(79,81){\circle*{5}} \put(81,84){1}
\put(93,54){\circle*{5}} \put(95,57){1}
\put(60,118){\line(0,1){17}}
\put(60,118){\circle*{5}} \put(62,116){3}
\put(60,138){\circle*{5}} \put(62,138){3}
\put(30,10){Figure 1a}
             
\put(127,54){\line(1,2){33}} \put(193,54){\line(-1,2){33}}
\put(127,54){\circle*{5}} \put(131,57){1} 
\put(141,81){\circle*{5}} \put(145,84){1}
\put(179,81){\circle*{5}} \put(181,84){1}
\put(193,54){\circle*{5}} \put(195,57){1}
\put(160,118){\line(0,1){17}}
\put(160,118){\circle*{5}} \put(162,116){2}
\put(160,138){\circle*{5}} \put(162,138){2}
\put(196,84){\line(-1,1){34}} \put(197,81){\circle*{5}}
\put(201,84){1}
\put(130,15){Figure 1b}
 \end{picture}\] 
\end{example}

Now we are in position to state and prove the main result of this paper.
\begin{theorem} \label{corr}
For any positive integer $n$ and for any positive integer $k\leq n$, there is
a bijection $B$ from the set 
of 2-stack sortable $n$-permutations with $k$  runs onto that
of \btrees on $n+1$ nodes having $k$ leaves.  Therefore,
 $W_{(n,k)}=T_{(n,k)}$. Moreover, for any 2-stack sortable $n$-permutation
$p$, the label of the root of $B(p)$ equals
$rl(p)$.  \end{theorem} 

\begin{proof} We are going to prove the statements of the Theorem by
induction on $n$ and $k$. For $k=1$ and any $n$, the statements are
clearly true: there is one $n$-permutation with one run, namely
$12\cdots n$, and there is one \btree having only one leaf, namely the
one consisting of a single path. The root label of this tree is 1, and
the number of \ril maxima of this permutation is 1, so the initial case
is proved.

Now let $k>1$ and suppose we already know the statements for all
permutations  shorter than $n$.   
Let $p$ be any 2-stack sortable $n$-permutation with $k$ runs.

\begin{itemize}
\item
 First suppose $p$ is of 
type 1. Consider $F(p)$, where $F$ is the bijection defined in the proof
of Lemma \ref{egyes}. Forget for now which  vertex of $F(p)$ is marked, then
$F(p)$ is a 2-stack sortable $(n-1)$-permutation with $k$ runs. So, by
our induction hypothesis, $B$ associates a \btree $T'$ on $n$ nodes to it,
which has $k$ leaves and whose root label is equal to $rl(F(p))\geq
rl(p)$. Now add one node above the root $x$ of $T'$ and ``recall'' which
\ril maximum of $F(p)$ was marked.  If it was the $t$th \ril maximum,
that is, if $p$ had exactly $t$ \ril maxima,
then let the label of this new root $x$ (and necessarily, that of the
old root, 
which is the only child of $x$) be $t$. Then we certainly get a \btree $T$
as $t\leq rl(F(p))$. We set $B(p)=T$. By the argument of Lemma
\ref{egyes} we see that $p$ can be recovered from $T$ (as $T'$, and thus
$F(p)$ can), and that $T$ has $k$ leaves, and its root label is $rl(p)$
as it should be. 
\item If $p=\epsilon$, then $B(p)$ is the only \btree on two nodes.
We point out that here, as well as in the previous case, the root of
$B(p)$ has only one child.

\item Now suppose $p$ is of type 2. Take $H(p)=(q_1,q_2,\cdots ,q_r)$ as
defined in Corollary \ref{szet} and consider $B(q_i)$ for all $i$ as
defined by the previous two cases. $B(q_i)$ has $q_i+1$ nodes, so by
contracting the roots of all $B(q_i)$ to one node (keeping the $B(q_i)$
ordered from left to right), we get a tree $T$ with $1+\sum_{i=1}^r
|q_i|=n+1$ nodes. We set $B(p)=T$. The components $q_i$ can be recovered
from $T$ by simply cutting off its root and adding a new root to the top
of each component obtained. (Recall that in all the $B(q_i)$, the root has
only one child). The number of leaves of $T$ equals the sum of the
leaves of the $B(q_i)$, so by the previous two cases and  by Corollary
\ref{szet},  it equals the number of runs of $p$. The root label of $T$
is by definition the sum of the labels of its children, which is,
again by the previous two cases and by Corollary \ref{szet}, equal to
$rl(p)$.  
\end{itemize}
This completes the proof of the Theorem.
\end{proof}

We remark that our algorithm certainly sets up a bijection between the
set of {\em all} 2-stack sortable $n$-permutations and that of {\em all}
\btrees on $n+1$ vertices. 

\section{New enumerative results}
Exercise 2.9.8.b in \cite{jackson} shows that the number of nonseparable
rooted planar maps with $f+1$ faces and $p+1$ vertices is equal to 
\begin{equation} \label{dmj} \frac{(2f+p-2)!(2p+f-2)!}{f!p!(2f-1)!(2p-1)!}.
\end{equation}

It is shown in \cite{cori} that the number of these maps equals the
number of \btrees on $n+1$ nodes with $k=f+1$ leaves and $n+1-k=p+1$
internal nodes. (Here the root is counted as an internal node). By
Theorem \ref{corr} this implies the result we announced in the
introduction. 

\begin{theorem} \label{enum} Let $W_{(n,k)}$ be the number of 2-stack sortable 
$n$-permutations  with $k$ runs. Then for all $1\leq k\leq n$ we have \[
W_{(n,k)}=\frac{(n+k-1)!(2n-k)!}{k!(n+1-k)!(2k-1)!(2n-2k+1)!}.\]
\end{theorem}

In particular, equation (\ref{dmj}) shows that the roles of $f+1=k$ and
$p+1=n+1-k$ are symmetric. In other words, $W_{(n,k)}=W_{(n,n+1-k)}$. A
permutation with $k$ runs has $k-1$ descents, so have proved the
following surprising Corollary.

\begin{corollary} The number of 2-stack sortable $n$-permutations with
$k$ descents equals that of those with $k$ ascents. \end{corollary}

A routine computation yields that $W_{(n,k)}/W_{(n,k-1)}>1$ if and only if
$k\leq (n+1)/2$, and we already know from the previous Corollary that if
$n=2k$, then $W_{(n,k)}=W_{(n,k+1)}$. So we have proved the following
Corollary. 

\begin{corollary} For fixed $n$, the sequence $W_{(n,k)}$ is unimodal in
$k$, and 
its peak is at $k=[(n+1)/2]$. \end{corollary}   
\section*{Acknowledgement}

I am grateful to Gilles Schaeffer
who sent me the preprint \cite{cori}.

\penalty-5000

\end{document}